\documentclass[12pt]{article}
\usepackage{graphicx}
\usepackage{amsmath}
\usepackage{amssymb}
\usepackage{theorem}

\numberwithin{equation}{section}

\newtheorem{thm}{Theorem}[section]

{\theorembodyfont{\rmfamily}

\newtheorem{example}[thm]{Example}

\newtheorem{rmk}[thm]{Remark}
}

\newcommand{\qed}{\hfill \mbox{\raggedright \rule{.07in}{.1in}}}

\newenvironment{proof}{\vspace{1ex}\noindent{\bf
Proof}\hspace{0.5em}}{\hfill\qed\vspace{1ex}}

\newcommand{\Eg}{E_{\Gamma}}
\newcommand{\R}{{\mathbb R}}

\newcommand{\N}{{\mathbb N}}

\newcommand{\SMALL}{\textstyle}

 \def\be{\begin{equation}}
 \def\ee{\end{equation}}
 
 \def\bea{\begin{eqnarray}}
 \def\eea{\end{eqnarray}}

\title{Decay of correlations and invariance principles \\ for dispersing
billiards with cusps, \\ and related planar billiard flows}

\author{
P\'eter B\'alint~\thanks{Institute of Mathematics,
Budapest University of Technology and Economics,
Budapest, Hungary, pet@math.bme.hu}
\and
Ian Melbourne~\thanks{Department of Mathematics,
University of Surrey, Guildford GU2 7XH, UK, ism@math.uh.edu}
}

\date{13 May 2008}

\begin{document}

\maketitle

\begin{abstract}
Following recent work of Chernov, Markarian, and Zhang, it is known
that the billiard map for dispersing billiards with zero angle
cusps has slow decay of correlations with rate $1/n$.
Since the collisions inside a cusp occur in quick succession, it is
reasonable to expect a much faster decay rate in continuous time.
In this paper we prove that the flow is {\em rapid mixing}:
correlations decay faster than any polynomial rate.
A consequence is that the flow admits strong statistical properties such as
the almost sure invariance principle, even though the billiard map
does not.

The techniques in this paper yield new results for other standard examples
in planar billiards, including Bunimovich flowers and stadia.
\end{abstract}

\section{Introduction}
\label{sec-intro}

Lorentz gas models (and the associated discrete time billiard maps) are
an important class of examples in mathematical physics.   Their
systematic study from the viewpoint of smooth ergodic theory was begun by
Sina{\u\i}~\cite{Sinai70} in the 1970's.
Sina{\u\i} focused on planar periodic Lorentz gases
with convex obstacles (dispersing billiards) proving properties such as
uniform hyperbolicity and ergodicity.   The situation is analogous to
the much-studied case of Anosov and Axiom~A systems for which advanced
statistical properties such as decay of correlations
and central limit theorems are by now classical~\cite{Bowen75,Ruelle78,Sinai72}.
In particular, if $\Lambda\subset M$ is a mixing
hyperbolic basic set with Gibbs measure $\mu$
for an Axiom~A diffeomorphism $f:M\to M$, then
for all H\"older observables $v,w:\Lambda\to\R$ the correlation function
\[
\textstyle \rho(n)=\int_\Lambda v\,w\circ f^n\,d\mu
-\int_\Lambda v\,d\mu \int_\Lambda w\,d\mu,
\]
decays exponentially: $\rho(n)=O(\tau^n)$ for some $\tau\in(0,1)$.
The corresponding results for billiards are complicated by the
presence of singularities and unbounded derivatives, but these difficulties
were eventually overcome by Young~\cite{Young98} and
Chernov~\cite{Chernov99} who proved exponential decay of correlations
for planar periodic dispersing billiards by constructing a certain type
of infinite Markov extension called a Young tower~\cite{Young98,Young99}.

The situation for flows is considerably more subtle.   Whereas mixing
Axiom~A diffeomorphisms always have exponential decay of correlations
for H\"older observables, there exist Axiom~A flows that are mixing but
at an arbitrarily slow rate~\cite{Ruelle83,Pollicott85}.   Currently, it is
not known if exponential decay of correlations is typical for Axiom~A flows.
However, Dolgopyat~\cite{Dolgopyat98b} (see also~\cite{FMT07})
proved that typical Axiom~A flows are rapid mixing (faster than any
polynomial rate) for sufficiently regular observables.
Melbourne~\cite{M07,Mapp} extended this result to flows with a
Poincar\'e map modelled by a Young tower, proving for example that
the planar periodic Lorentz gas with finite horizon is rapid mixing.
For this particular flow, Chernov~\cite{Chernov07} was able to prove
stretched exponential decay of correlations.

Thus, generally speaking
our understanding of planar periodic dispersing billiards and Lorentz gases is now as advanced as in the Axiom~A case.
However, there are a number of classes of planar billiards for which
open questions remain~\cite{ChernovDolgopyat06,ChernovMarkarian,ChernovYoung00}.    In this paper, we focus mainly on questions related to
decay of correlations for Lorentz gas flows, building on~\cite{M07}
(rapid mixing) and~\cite{Mapp} (slow mixing).

\begin{example}[Dispersing billiards with cusps]
\label{ex-cusps}
Dispersing billiards
with cusps, where the boundary curves are all dispersing
but the interior angles at corner points are zero, turned out to be
much more involved technically than the usual Lorentz gases. For example
ergodicity for the dispersing billiard map with cusps was only proved in the mid
1990's by Reh\'a{\v c}ek~\cite{Rehacek}. By standard arguments, ergodicity implies K-mixing~\cite[Chapter~6]{ChernovMarkarian} and the Bernoulli property
\cite{ChernovHaskell,OrnsteinWeiss} for all the hyperbolic billiard examples.

Concerning rates of mixing, Chernov \& Markarian~\cite{ChernovMarkarian07}
proved, using the method of Young towers, that decay of correlations for the
billiard map satisfies $\rho(n)=O((\log n)^2/n)$.  It was anticipated that
the logarithmic factor is an artifact of the proof, and this factor was removed
by Chernov \& Zhang~\cite{ChernovZhang08} yielding the decay rate
$\rho(n)=O(1/n)$.   It is believed that this result is optimal for the billiard
map.
This rate of decay of correlations is too weak for strong statistical limit laws
(see~\cite{BalintGouezel06}).

Since the collisions inside a cusp occur in quick succession, it is
reasonable to expect a much faster decay rate in continuous
time~\cite[Section~5.6]{ChernovYoung00}, possibly even exponential decay of
correlations.   Exponential decay seems to be beyond current techniques even for
simpler problems (such as the finite horizon planar periodic Lorentz gas).
Previously, no results on decay of correlations for the flow
(not even slow mixing) were available when there are cusps.
In this paper \textit{we prove rapid mixing for the flow}.

 A byproduct of our proof is the {\em almost sure invariance
principle (ASIP) for the flow}.
We note that the ASIP implies numerous statistical limit laws such as the central
limit theorem, the law of the iterated logarithm, and their functional versions,
see for example~\cite{PhilippStout75,MNapp}.
\end{example}

\begin{example}[Bunimovich flowers]
\label{ex-flowers}
It was discovered by Bunimovich in the 1970's that billiard tables with
focusing boundary components may also show hyperbolic behaviour. For
the first  examples of such tables, constructed for example in \cite{Bunimovich74},
the boundary components are either dispersing, or focusing arcs of circles, subject to
some further constraints (see Section~\ref{subsec-flowers} below). Given their typical
shape, such billiards are often called Bunimovich flowers. Later the probably
best known example of a hyperbolic billiard with focusing boundary components, the
Bunimovich stadium, was introduced in~\cite{Bunimovich79}. Focusing hyperbolic
billiards have a rich literature by now -- it is worth pointing out Wojtkowski's
contribution, \cite{Wojtkowski} in particular, which clarified much of the mechanism
behind hyperbolic behaviour. For further historical accounts, and detailed proofs of
ergodicity, we refer to the monograph \cite{ChernovMarkarian} and references therein.

As for rates of mixing, Chernov \& Zhang~\cite{ChernovZhang05} show that in Bunimovich flowers
the billiard map has decay of correlations $\rho(n)=O((\log n)^3/n^2)$.
(Again, the logarithmic factor appears to be an artifact of the proof and it is
expected that $1/n^2$ is the optimal rate.)
It follows that the map satisfies statistical limit laws such as the ASIP.
See~\cite{MN05} for the scalar case and~\cite{MNapp} for $\R^d$-valued observables.

By~\cite{MT04}, the billiard flow
immediately inherits the ASIP (and hence the other limit laws).
However, previous results were unable to establish estimates for
decay of correlations of the flow.
We recall from \cite{ChernovZhang05} that the only effect slowing down the
decay rate of the billiard map is sliding along a circular arc, where collisions
occur in quick succession.
Just as in the cusp example, it is reasonable to expect that
the flow mixes better than the map in Bunimovich flowers.

In this paper, {\em we prove rapid mixing for the flow}.
\end{example}

\begin{example}[Bunimovich stadia]
\label{ex-stadia}
Bunimovich~\cite{Bunimovich79} established hyperbolicity and ergodicity
for the stadium billiard bounded by two parallel lines connecting two semicircles.
Markarian~\cite{Markarian04} proved decay of correlations
$O((\log n)^2/n)$ for the billiard map, and Chernov \& Zhang~\cite{ChernovZhang08}
improved this to obtain the optimal rate $1/n$.
This is too weak for strong statistical limit laws;
indeed B\'alint \& Gou\"ezel~\cite{BalintGouezel06} prove a nonstandard limit
law (nonstandard domain of attraction of the normal distribution with
$\sqrt{n\log n}$ normalization) for typical observables.
The same limit law holds for the flow~\cite{BalintGouezel06,MT04}.
In particular, the ASIP fails for both the billiard map and the flow.

This time, we do not expect the flow to mix more rapidly than the map,
but slow mixing does not follow from previous results.
In this paper, {\em we prove that the stadium flow decays at the same
slow rate $1/t$ as the map}.
\end{example}

\begin{rmk}
Dolgopyat's results~\cite{Dolgopyat98b} on rapid mixing for Axiom~A flows
require the observables to be smooth in the flow direction.
This restriction is inherited by the generalisation in~\cite{M07,Mapp} to
nonuniformly hyperbolic flows, and by the current paper.
In particular, our results on rates of mixing for Lorentz gases do not apply
to certain physically relevant observables such as position and velocity.
At present, the only result on decay of correlations for Lorentz gas flows
that includes general H\"older observables is Chernov's result
for the planar periodic Lorentz gas with finite horizon~\cite{Chernov07}.

On the other hand, our proof of the ASIP for billiard flows with cusps
in Example~\ref{ex-cusps} is valid for all H\"older observables.
\end{rmk}

The remainder of the paper is organised as follows.
In Section~\ref{sec-Young}, we recall background material on Young towers~\cite{Young98,Young99}, and on decay of correlations for flows possessing a Poincar\'e map
modelled by a Young tower~\cite{M07,Mapp}.
In Section~\ref{sec-bill}, we present our new results for
Examples~\ref{ex-cusps},~\ref{ex-flowers} and~\ref{ex-stadia}.

\section{Young towers and flows}
\label{sec-Young}

Let $f:M \to M$ be a map, and let $h:M \to \R^+$ be a roof function.
Let $\phi_t:M^h\to M^h$ be the corresponding suspension flow.
Thus the map $f(x)=\phi_{h(x)}(x)$ is a Poincar\'e map for the flow.
Conversely, given $f$ and $h$ we define the suspension
$M^h=\{(x,u)\in M\times\R:0\le u\le h(x)\}/\sim$
where $(x,h(x))\sim(fx,0)$.  Then the suspension flow is given by
$\phi_t(x,u)=(x,u+t)$ computed modulo identifications.

In our examples, $M$ is a finite union of Riemannian manifolds.
We will say that the map
$f:M\to M$ is \textit{modelled by a Young tower} if it satisfies
the axioms introduced by Young~\cite{Young98,Young99}. In particular,
there is a set $Y\subset M$ with ${\rm Leb}(Y)>0$ that possesses an appropriate
hyperbolic product structure. Furthermore there exists a measurable
partition of $Y$ into countably many sets $Y_i$, and
$r:Y\to \N$ with constant values $r_{Y_i}$ on the $Y_i$, such that the induced
map $F=f^r:Y\to Y$ is smooth and uniformly hyperbolic when restricted to any
$Y_i$, with respect to the hyperbolic product structure of $Y$. It is important to note that $r$ is
not necessarily the first return time to $Y$, thus the corresponding tower map is a Markov extension
of the original map $f:M\to M$. For further details of the construction
see the original references.

For examples from billiards, there is a natural invariant measure $\mu$ (the Liouville
measure) that is equivalent to Lebesgue measure on $M$ with $L^\infty$ density.
In what follows we always assume that both the discrete return time $r$ (defined on $Y$) and the
roof function $h$ (defined on $M$) are integrable with respect to Lebesgue measure
and hence $\mu$.
The natural flow-invariant probability measure $\mu^h$ on $M^h$ given by
$\mu^h=\mu\times{\rm lebesgue}/\int_M h\,d\mu$ coincides with Lebesgue measure
for the billiard examples.

We aim to compute the rate of decay of correlations and statistical limit
laws for sufficiently regular observables $v,w:M^h\to\R$.
For $m\ge0$, $\eta\in(0,1)$, the function space $C^{m,\eta}(M^h)$ consists of those $v:M^h\to \R$
for which $\|v\|_{m,\eta}=\|v\|_{\eta}+\|\partial_t v\|_{\eta}+\dots+\|\partial_t^m v\|_{\eta}
<\infty$ where $\partial_t$ denotes the derivative in the flow direction and
$\|v\|_{\eta}=|v|_\infty + \sup_{x\neq y,\, u\in[0,\max\{h(x),h(y)\}]}|v(x,u)-v(y,u)|/d(x,y)^\eta$.

We now summarize results from~\cite{M07,Mapp} used in Section~\ref{sec-bill}.  Define
\[
\textstyle \rho_{v,w}(t)=\int_{M^h} v\,w\circ \phi_t\,d\mu^h
-\int_{M^h} v\,d\mu^h \int_{M^h} w\,d\mu^h.
\]

\begin{thm}[ \cite{M07} ]
\label{thm-rapid}
Let $f:M \to M$ be a map modelled by a Young tower and
$h:M \to \R^+$ be a H\"older roof function.   Assume
\begin{itemize}
\item Exponential tails: ${\rm Leb}(y\in Y:r(y)>n)=O(\gamma^n)$ for some $\gamma\in(0,1)$.
\end{itemize}

Then typically the flow
 $\phi_t:M^h\to M^h$ is {\em rapid mixing}: there exists $\eta\in(0,1)$ and
for any $n\ge1$, there exists $m\ge0$, $C\ge1$, such that
\[
|\rho_{v,w}(t)|\le C\|v\|_{m,\eta}\|w\|_{m,\eta}\, t^{-n},
\]
for all $v,w\in C^{m,\eta}(M^h)$ and all $t>0$.
\end{thm}

\begin{thm}[ \cite{Mapp} ]
\label{thm-slow}
Let $f:M \to M$ be a map modelled by a Young tower and
$h:M \to \R^+$ be a H\"older roof function.   Assume
\begin{itemize}
\item Polynomial tails:
${\rm Leb}(y\in Y:r(y)>n)=O(1/n^{\beta+1})$ for some $\beta>0$, and
\item $1/h\in L^\infty$, or more generally there
is an integer $p\ge1$ such that $1/h_p\in L^\infty$
where $h_p=h+h\circ f+\dots+h\circ f^{p-1}$.
\end{itemize}

Then typically the flow
 $\phi_t:M^h\to M^h$ is {\em polynomial mixing with rate $1/t^{\beta}$}:
there exists $\eta\in(0,1)$, $m\ge0$, $C\ge1$, such that
\[
|\rho_{v,w}(t)|\le C\|v\|_{m,\eta}\|w\|_{m,\eta}\, t^{-\beta},
\]
for all $v,w\in C^{m,\eta}(M^h)$ and all $t>0$.
\end{thm}

\begin{rmk}  In general, Theorems~\ref{thm-rapid} and~\ref{thm-slow} hold
{\em typically} subject to a nondegeneracy condition.
However, this nondegeneracy condition
is automatic for the billiards examples considered here
because of the contact-like structure~\cite{Mapp}.
\end{rmk}

\begin{rmk}   Theorem~\ref{thm-rapid} implies that
in situations where Young~\cite{Young98} obtains exponential decay
of correlations for the billiard map $f:M\to M$,
rapid decay of correlations holds for the flow $\phi_t$.   Similarly,
Theorem~\ref{thm-slow} implies that
in situations where Young~\cite{Young99} obtains polynomial decay
of correlations for the billiard
map $f:M\to M$ with rate $O(1/n^\beta)$,
polynomial decay of correlations with rate $O(1/t^\beta)$ holds for the flow
$\phi_t$.
\end{rmk}

\begin{rmk}
\label{rmk-rapid}
The condition that $h$ is H\"older can be relaxed in
Theorems~\ref{thm-rapid} and~\ref{thm-slow}.
It suffices that $h$ is \textit{uniformly piecewise H\"older} as follows.
The billiard examples are examples of smooth dynamical systems with
at most \textit{countably many singularities}.  Write
$M=\bigcup M_j$ where the $M_j$ are the maximal subsets of $M$ on which
$f$ is smooth.  It is enough that $\sup_{j}\|h1_{M_j}\|_\eta<\infty$.

By construction, the Young tower that models the billiard map ``respects'' the partition
$\bigcup M_j$ in the following sense: as $F\vert_{Y_i}$ is smooth, $\bigcup Y_i$
is necessarily a refinement of the restriction of $\bigcup M_j$ onto $Y$.
\end{rmk}

Let us make one more remark. Even though this is a very simple observation, this is the
\textit{Key Observation} behind the arguments of Section~\ref{sec-bill}.
Note that $M$ can be viewed as a Poincar\'e cross-section to the flow.
It may be the case that there is an alternative cross-section $\widehat M$
with Poincar\'e map $\widehat f:\widehat M \to \widehat M$ (and roof function $\widehat h:\widehat M\to\R^+$)
that can also be modelled by a Young tower and such that
the tail distribution
of the Young tower for $\widehat f:\widehat M \to \widehat M$
allows for a better control of the statistical properties.
Note that the corresponding suspension flow $\widehat\phi_t:\widehat M^{\widehat h}\to \widehat M^{\widehat h}$
is identical to the original flow.

\subsubsection*{Statistical limit laws}

\begin{thm}[ \cite{MNapp} ]
\label{thm-ASIP}
Let $f:M \to M$ be a map modelled by a Young tower and
$h:M \to \R^+$ be a H\"older roof function.   Assume
\begin{itemize}
\item Summable tails:
${\rm Leb}(y\in Y:r(y)>n)=O(1/n^{\beta+1})$ for some $\beta>1$.
\end{itemize}

Then the {\emph (vector-valued) ASIP holds}: there exists $\lambda>0$, and for any
H\"older observable $v:M^h \to \R^d$ with mean zero
there is a $d$-dimensional Brownian motion $W$, such that
(on a possibly enriched probability space)
\[
{\SMALL\int}_0^T v\circ\phi_t\, dt = W(T) + O(T^{\frac12-\lambda}) \enspace
\text{as $T\to\infty\enspace$ a.e.}
\]
\end{thm}

\begin{rmk}
Again, it suffices that the roof function $h$ is uniformly piecewise H\"older
as in Remark~\ref{rmk-rapid}.
Similarly, the H\"older regularity requirement for the observable $v$ in
Theorem~\ref{thm-ASIP} can be relaxed.
Define $V(x)=\int_0^{h(x)}v\circ\phi_t\,dt$.  It suffices that
$V: M\to\R^d$ is uniformly piecewise H\"older.
\end{rmk}

\begin{rmk}   \label{rmk-time1}
A more difficult question is the ASIP for the time-one map of a flow.
A general argument of~\cite{FMT03,MT02} yields the scalar ASIP (for real-valued
observables) in situations where (i) summable decay of correlations is established
by certain techniques, (ii) the class of dynamical systems is closed under
time reversal.   This includes time-one maps for Axiom~A flows~\cite{MT02}.
Again the result is restricted to observables that are sufficiently smooth in the
flow direction.

It is immediate from these considerations that the scalar ASIP holds for time-one maps
of the flows in the cusp and flower examples (Examples~\ref{ex-cusps} and~\ref{ex-flowers}).   The vector-valued ASIP remains an open question even for time-one maps of
Axiom~A flows.
\end{rmk}

\section{Billiard flows}
\label{sec-bill}

In this section, we specialise to billiard flows and treat the three examples
discussed in the introduction.   We begin with Example~\ref{ex-flowers} (flowers) since
this is the simplest.  Then we treat Examples~\ref{ex-cusps} (cusps)
and~\ref{ex-stadia} (stadia).

Let $Q$ be a compact domain in $\R^2$ with piecewise $C^3$ boundary $\partial Q$
and let us consider the corresponding billiard dynamics in $Q$ -- the
motion of a point particle that travels with unit speed and bounces off
$\partial Q$ elastically (the angle of incidence is equal to the angle of reflection).
The resulting flow $\phi_t$ is a three-dimensional volume preserving
flow.   There is a natural two-dimensional cross-section
$M=\partial Q\times[-\pi/2,\pi/2]$ corresponding to collisions, and the
Poincar\'e map $f:M\to M$ is called the billiard map.
We choose coordinates $x=(r,\varphi)$ on $M$, where $r$ (the position)
is the arclength parameter along $\partial Q$ and $\varphi$ (the outgoing velocity)
is the angle made by the reflected trajectory and the normal to the boundary.
The natural invariant measure for the billiard map is is the Liouville measure
$d\mu=\cos\varphi\, dr\,d\varphi$.

By components or arcs we mean the maximal connected $C^3$-smooth pieces of the boundary $\partial Q$.
We use the same terminology for the corresponding (two-dimensional)
components of the phase space $M$.
We assume that the curvature of any component has a fixed sign, thus we may consider convex
inwards (dispersing), convex outwards (focusing) and neutral arcs, and decompose the phase space as
$M=M^+\cup M^-\cup M^0$, respectively.

The dynamics of the billiard map/flow may depend sensitively on further geometrical properties of
$\partial Q$, including the types of the components present and their distances and angles.
Below we restrict to certain physically relevant examples and cite the literature for their
known dynamical properties. For further details on billiard dynamics in general, see \cite{ChernovMarkarian} and
references therein.

Our analysis relies on the standard
fact that the billiard map $f$ and the roof function
(collision time) $h$ are piecewise H\"older continuous.   A simple argument shows that
the H\"older exponent is $\frac12$ for both $f$ and $h$.
Note that $h(x)=|f(x)-x|$
so it suffices to verify that $f$ is $\frac12$-H\"older.

It is clear that $f$ is smooth except near preimages of tangencies
(grazing collisions) and writing $f=(f_r,f_\varphi)$ a standard
calculation~\cite[Section~2.11]{ChernovMarkarian} shows that $f'$ is bounded
except for singularities where it behaves like $1/\cos(f_\varphi(x))$.
Writing $\psi=\pi/2-\varphi$ and suppressing the $r$-coordinates,
we have $f'_{\psi}(x)\sim 1/\sin(f_{\psi}(x))\sim 1/(f_{\psi}(x))$
so that $((f_{\psi}(x))^2)'\sim 1$. Hence $f_{\psi}(x)\sim\psi^{1/2}$
which is $\frac12$-H\"older as required.

\subsection{Bunimovich flowers}
\label{subsec-flowers}
Bunimovich~\cite{Bunimovich74} studied hyperbolicity and ergodicity for a class of billiard tables $Q\subset\R^2$ such that (i) $M^0=\emptyset$ (no neutral components); the dispersing components $M^+$
may have arbitrary geometry, in contrast (ii) $M^-$ consists of circular arcs all strictly smaller
than a semicircle; (iii) if such a circular arc was extended to a full circle, this circle would be contained
in $Q$; (iv) the neighbouring components are transversal (no cusps). Actually, the conditions can be relaxed
and even neutral components can be allowed as long as certain technical assumptions are satisfied, see
\cite{ChernovZhang05} for details.

Chernov \& Zhang~\cite{ChernovZhang05} show that for Bunimovich flowers the billiard map
$f:M \to M$ can be modelled
by a Young tower with tail estimate roughly $O(1/n^3)$.
It then follows from Young~\cite{Young99} that the
map has decay of correlations roughly $O(1/n^2)$.
By~\cite{MNapp}, the vector-valued ASIP holds for H\"older observables
for both the billiard map and the billiard flow.

However, there is no immediate conclusion for the
decay of correlations of the
flow.  (Exponential tails are required in Theorem~\ref{thm-rapid}
while Theorem~\ref{thm-slow} would yield at best
$O(1/t^2)$ decay.    In any case, the assumption that
$h$ is bounded below in Theorem~\ref{thm-slow} is violated.)
Nevertheless, we now prove rapid mixing for the flow for sufficiently smooth
observables.
(The scalar ASIP for the time-one map for the flow then follows from
Remark~\ref{rmk-time1}.)

\begin{thm} \label{thm-flowers}
Let $\phi_t$ be the flow corresponding to a Bunimovich flower.
Then $\phi_t$ is rapid mixing (in the sense of Theorem~\ref{thm-rapid}).
\end{thm}

\begin{proof}
The argument in
Chernov \& Zhang~\cite[p.~1546]{ChernovZhang05} demonstrates the existence
of an alternative cross-section $\widehat M\subset M$ (with
Poincar\'e map $\widehat f:\widehat M \to \widehat M$) that is modelled by a Young tower
with exponential tails.   We claim that the corresponding collision time
$\widehat h:\widehat M \to \R^+$ satisfies the uniformly piecewise
H\"older condition in Remark~\ref{rmk-rapid},
expressed of course in terms of $\widehat M,\widehat f,\widehat h$.
Rapid mixing follows from Theorem~\ref{thm-rapid}.

First, we recall the definition of $\widehat M$ in~\cite{ChernovZhang05}.
As mentioned above we write $M=M^+\cup M^-$ where
 $M^+$ corresponds to dispersing arcs and
 $M^-$ corresponds to focusing arcs.   Then $\widehat M$ is given by $\widehat M=M^+\cup E$
 where $E\subset M^-$ consists of \textit{only the first (sliding) collisions at each
 focusing arc}, so $E=M^-\cap f(M^+)$.

  Next we verify that $\widehat h\in L^\infty$.
Since $\widehat h=h$ on dispersing arcs we can restrict attention to a single
focusing arc $\Gamma(\subset M^-)$, and in particular the single first collision
set $\Eg =\Gamma\cap E$.  Consider $x\in \Eg$.
 There is an integer $n\ge2$ such that $f^jx\in\Gamma$
 for $j=1,\dots,n$ and $f^{n+1}x\not\in\Gamma$.  Hence
\begin{align}
\label{def:g}
\widehat h=g+h\circ f^n,
\end{align}
where $g(x)=h(x)+h(fx)+\dots+h(f^{n-1}x)$ is the amount of time it takes
to ``slide'' along the arc $\Gamma$.
Hence $|\widehat h1_E|_\infty\le |\Gamma|+|h|_\infty$ and so  $\widehat h$ is bounded.

Finally, we verify that there is a uniform H\"older constant
for $\widehat h$ on each partition element. Again we may restrict to a single first collision set
$\Eg\subset\Gamma$ where $\Gamma$ is a focusing arc, but this time
it is necessary to consider the finer partition $\Eg=\bigcup_{n\ge1} E_n$ where
$E_n=\{x\in E:f^jx\in\Gamma\enspace\text{for $j=1,\dots,n$ and}\,
f^{n+1}x\not\in\Gamma\}$.   (In other words, $\widehat f=f^{n+1}$ on $E_n$.)
Note that the partition elements $E_n$ coincide with the smoothness
components of $\widehat f:\widehat M\to \widehat M$, which is in agreement with
Remark~\ref{rmk-rapid}.

Let $R$ denote the radius of the focusing arc $\Gamma$.
Let $(r,\varphi)$ denote the standard coordinates and note that
sliding occurs for $\varphi\sim\pi/2$ so it is convenient to introduce
$\psi=\pi/2-\varphi$.
It is immediate
from the geometry of the circle that
\[
f(x)=(r+2R\psi,\psi), \quad h(x)=R\sin\psi
\]
for all $x=(r,\psi)\in \Gamma\cap f^{-1}\Gamma$.
Hence
\begin{align} \label{eq-flower_est1}
|g(x_1)-g(x_2)|\le nR|\sin\psi_1-\sin\psi_2| \le nR|\psi_1-\psi_2|
\end{align}
for all $x_1,x_2\in E_n$.
On the other hand, since the angle $\psi$ remains constant during the $n$ grazing iterations, it
is elementary that there is a constant $C_1$ depending only on
$\Gamma$ such that $C_1/(n+1)\le \psi\le C_1/n$ for all
$x\in E_n$.
Hence $|\psi_1-\psi_2|\le C_1(\frac1n-\frac{1}{n+1})\le C_1/n^2$ and so
\begin{align} \label{eq-flower_est2}
|\psi_1-\psi_2|^\frac12 \le C_1^\frac12 n^{-1}.
\end{align}
Combining estimates~\eqref{eq-flower_est1} and~\eqref{eq-flower_est2}, we
obtain $|g(x_1)-g(x_2)|\le C_2|\psi_1-\psi_2|^\frac12$ where
$C_2=RC_1^\frac12$ depends only on the arc $\Gamma$.

In addition, since $f^n(r,\psi)=(r+2nR\psi,\psi)$ on $E_n$ and
$h$ $\frac12$-H\"older, $|h(f^nx_1)-h(f^nx_2)|\le C_3(|r_1-r_2|+(2Rn+1)|\psi_1-\psi_2|)^\frac12
\le C_4(|r_1-r_2|+|\psi_1-\psi_2|^\frac12)^\frac12
\le C_4(|r_1-r_2|^\frac12+|\psi_1-\psi_2|^\frac14)$.
Combining the estimates for $g$ and $h\circ f^n$,  we have shown that
$\widehat h$ is H\"older on $E_n$ with constant independent of $n$ (and exponent $\frac14$).
Since there are only finitely many arcs, this completes the proof.
\end{proof}

\subsection{Dispersing billiards with cusps}

Chernov \& Markarian~\cite{ChernovMarkarian07} studied dispersing billiards
with cusps, where the boundary curves are all dispersing -- that is, $M=M^+$ --
but the interior angles at corner points are zero.
By~\cite{ChernovMarkarian07,ChernovZhang08}, the billiard
map $f:M\to M$ can be modelled
by a Young tower with tail estimate $O(1/n^2)$.
It follows from Young~\cite{Young99} that the
map has decay of correlations $O(1/n)$.
This is too weak for strong statistical limit laws (see~\cite{BalintGouezel06}).
Nevertheless we now prove rapid mixing and the ASIP for the flow.

\begin{thm} \label{thm-cusps}
Let $\phi_t$ be the flow corresponding to a billiard table with cusps.
Then $\phi_t$ is rapid mixing (in the sense of Theorem~\ref{thm-rapid})
and satisfies the ASIP (in the sense of Theorem~\ref{thm-ASIP}).
\end{thm}

\begin{proof}
Following Chernov \& Markarian~\cite{ChernovMarkarian07}
we define $\widehat M$ \textit{by excluding a neighbourhood of each cusp}.
The new collision map $\widehat f:\widehat M\to \widehat M$ is modelled by a Young tower with exponential
tails~\cite{ChernovMarkarian07}.
By Theorems~\ref{thm-rapid} and~\ref{thm-ASIP}, it again suffices to show that the new collision time $\widehat h:\widehat M \to \R^+$
satisfies the uniformly piecewise H\"older condition in Remark~\ref{rmk-rapid}.

Consider a single cusp and let $E_n$ be the set of those points in $\widehat M$ that spend
exactly $n$ iterates in the cusp before returning to $\widehat M$. Note again that, in accordance
with Remark~\ref{rmk-rapid}, the sets
$E_n$ coincide with the smoothness components of $\widehat f$.
The calculation in~\cite[p.~738]{ChernovMarkarian07}
shows that $\widehat h$ is bounded on $E_n$ independent of $n$.   Hence $\widehat h$ is
bounded on the whole of $\widehat M$.

It remains to find a uniform H\"older constant for $\widehat h$ on each $E_n$.
Explicit calculations are more difficult than in Section~\ref{subsec-flowers}
so we search for coarser estimates.
We claim that there are positive constants $\alpha_1$, $\alpha_2$, $\alpha_3$,
$C_1$, $C_2$ (independent of $n$) such that
\begin{itemize}
\item[(i)]
$|\widehat h(x_1)-\widehat h(x_2)|\le C_1n^{\alpha_1}|x_1-x_2|^{\alpha_2}$, and
\item[(ii)]
$|\widehat h(x_1)-\widehat h(x_2)|\le C_2n^{-\alpha_3}$,
\end{itemize}
for all $x_1,x_2\in E_n$.  It then follows that
$|\widehat h(x_1)-\widehat h(x_2)|^{\gamma+1}\le C|x_1-x_2|^{\alpha_2\gamma}$
with $\gamma=\alpha_3/\alpha_1$, $C=C_1^\gamma C_2$,
and so $\widehat h$ is uniformly piecewise H\"older as required, with exponent
$\alpha_2\gamma/(\gamma+1)$.

We verify the claim with $\alpha_1=\frac{11}{6}$, $\alpha_2=\frac14$
and $\alpha_3=\frac13$ (and hence H\"older exponent $\frac{1}{26}$).
Throughout, $C$ is a uniform constant that can change from line to line.
Let $x_1,x_2\in E_n$.  Recall that $f$ and $h$ are H\"older with exponent $\frac12$.
Furthermore, by~\cite[Proposition~4.1]{ChernovMarkarian07},
\[
|f^jx_1-f^jx_2|\le Cn^{5/3}|fx_1-fx_2|\le Cn^{5/3}|x_1-x_2|^{1/2}
\]
for $j=1,\dots,n$.
We deduce that
\begin{align*}
|\widehat h(x_1)-\widehat h(x_2)| &\le |h(x_1)-h(x_2)|+ \sum\limits_{j=1}^n |h(f^jx_1)-h(f^jx_2)| \\ 
& \le C|x_1-x_2|^{1/2}+ C\sum\limits_{j=1}^n |f^jx_1-f^jx_2|^{1/2}  \\ 
& \le C|x_1-x_2|^{1/2} + n Cn^{5/6}|x_1-x_2|^{1/4}\le Cn^{11/6}|x_1-x_2|^{1/4}
\end{align*}
establishing estimate (i).

By~\cite[p.~738, last line]{ChernovMarkarian07},
we have $h(fx_j)+\dots+ h(f^{n-1}x_j)\le Cn^{-1}$ for $j=1,2$.
As shown in~\cite[pp.~748--749]{ChernovMarkarian07}, the cell $E_n$
has diameter of order $n^{-2/3}$ so that $|h(x_1)-h(x_2)|\le C|x_1-x_2|^{1/2}
\le Cn^{-1/3}$.
By time-reversibility of the construction, $|h(f^nx_1)-h(f^nx_2)|\le Cn^{-1/3}$.
Altogether, we have $|\widehat h(x_1)-\widehat h(x_2)|\le Cn^{-1/3}$ establishing estimate (ii).
(The estimates of~\cite{ChernovMarkarian07} that we have used are established in
the first instance for a special billiard table with three cusps, but then
extended to the general situation in~\cite[Section~6]{ChernovMarkarian07}.)
\end{proof}

\begin{rmk}   Again, we stress that the vector-valued ASIP for the flow holds
for general (piecewise) H\"older observables, whereas rapid mixing is restricted
to sufficiently smooth observables (as is the scalar ASIP for the time-one
map of the flow which holds by Remark~\ref{rmk-time1}).
\end{rmk}
\subsection{Bunimovich stadia}
Bunimovich~\cite{Bunimovich79} established hyperbolicity and ergodicity
for the stadium billiard bounded by two parallel lines ($M^0$) connecting
two semicircles ($M^-$).
By~\cite{Markarian04,ChernovZhang08}, the billiard map
$f:M \to M$ can be modelled by a Young tower with tail estimate $O(1/n^2)$.
It follows from Young~\cite{Young99} that the
map has decay of correlations $O(1/n)$ for H\"older observables.
This is too weak for strong statistical limit laws;
indeed B\'alint \& Gou\"ezel~\cite{BalintGouezel06} prove a nonstandard limit
law (nonstandard domain of attraction of the normal distribution with
$\sqrt{n\log n}$ normalization) for typical observables.
By~\cite{MT04}, the same limit law holds for the flow.  In particular,
the ASIP fails for both the billiard map and the flow.

This time, we do not expect the flow to mix more rapidly than the map.
Nor can we apply Theorem~\ref{thm-slow} directly since the
roof function is not bounded below.  Nevertheless the conclusion
of Theorem~\ref{thm-slow} is valid.

\begin{thm} \label{thm-stadia}
Let $\phi_t$ be the flow corresponding to a Bunimovich stadium.
Then $\phi_t$ is polynomially mixing with rate $1/t$.
\end{thm}

\begin{proof}
The argument in
Chernov \& Zhang~\cite[p.~1548]{ChernovZhang05} demonstrates the existence
of an alternative cross-section $\widehat M\subset M$ (with
Poincar\'e map $\widehat f:\widehat M \to \widehat M$) that is modelled by a Young tower
with exponential tails.   However, the corresponding roof function
$\widehat h:\widehat M\to\R^+$ is unbounded, so
Theorem~\ref{thm-rapid} does not apply this time.
As mentioned above, they show that the corresponding tower for $f:M \to M$
has tails decaying at the rate $O(1/n^2)$.
We note that the tower for $M$ is built over the same base $Y$
as the tower for $\widehat M$ but is strictly higher (since it incorporates the returns
from $M$ to $\widehat M$).

Whereas the cross-section $\widehat M$ in~\cite{ChernovZhang05} takes account
of both sliding along circular arcs in  $M^-$ and bouncing between the parallel
straight lines in $M^0$, we define an intermediate cross-section
$\widehat M'$ with $\widehat M\subset \widehat M'\subset M$ that takes only account of sliding.
That is, we define $\widehat M'=M^0\cup E$ where $E=M^-\cap f(M^0)$ consists of the first (sliding)
collisions in $M^-$.

It is immediate that the tower for $\widehat f':\widehat M' \to \widehat M'$ has tail decay rate
no worse than the tower for $\widehat M$: it shares the same base $Y$
as the other two towers, but is lower than the tower for $\widehat M$
(and higher than the tower for  $M$).
In particular the tower for $\widehat M'$ has tail decay rate $O(1/n^2)$.

Let $\widehat h':\widehat M'\to\R^+$ denote the new roof function.
It remains to verify the hypotheses on $\widehat h'$ in Theorem~\ref{thm-slow}.
The proof that $\widehat h'$ is uniformly piecewise
H\"older is identical to that in Section~\ref{subsec-flowers}.
We claim that $\widehat h'_p$ is bounded below for $p=2$.
Indeed, $\widehat h'$ can approach zero only at points $x\in M^0$ which are close
to one of the endpoints of a straight boundary component, that is, to the corner
made with $M^-$, and have velocity almost tangent to the straight boundary
component.   But then $fx$ is about to undergo
a long sequence of sliding collisions so that $\widehat h'(fx)$ has magnitude
approximately the length of a semi-circular arc in $M^-$.
Hence $\widehat h'_2$ is bounded away from zero as required.~
\end{proof}

\paragraph{Acknowledgements}
This research was supported in part by EPSRC Grant EP/D055520/1; by OTKA (Hungarian
National Research Fund) grants: F 60206, TS 49835 and K 71693; and by the 
Bolyai scholarship of the Hungarian Academy of Sciences.
We would like to thank the Budapest University of Technology and Economics
 and the University of Surrey for
hospitality while part of this research was done.
IM is grateful to the University of Houston for the use of e-mail.

\end{document}